\documentclass[12pt]{amsart}
\usepackage{graphicx}
\usepackage{amssymb}
\usepackage{amsmath}

\newcommand{\C}{\mathbb C}

\newcommand{\D}{\mathbb D}
\newcommand{\F}{\mathbb Z}

\renewcommand{\P}{\mathbb P}

\newcommand{\Z}{\mathbb Z}

\newcommand{\cF}{\mathcal F}
\newcommand{\cI}{\mathcal I}
\newcommand{\cO}{\mathcal O}

\newcommand{\pa}{\partial}

\newcommand{\na}{\nabla}
\renewcommand{\a}{\alpha}

\newcommand{\De}{\mathit{\Delta}}
\newcommand{\lm}{\lambda}

\newcommand{\f}{\varphi}

\newcommand{\G}{\mathit{\Gamma}}

\newcommand{\W}{\mathit{\Omega}}
\newcommand{\w}{\omega}

\newcommand{\fX}{\frak{X}}

\newcommand{\tr}{\;^t}

\newcommand{\cX}{\widetilde X}
\newcommand{\wt}{\widetilde}
\newcommand{\wh}{\widehat}
\newcommand{\bu}{\bullet}

\newcommand{\cE}{\mathcal{E}}

\newcommand{\pr}{\frak{p}}
\newcommand{\pro}{\mathit{pr}}

\newcommand{\Qed}{\qed\medskip}
\newtheorem{theorem}{Theorem}[section]

\newtheorem{proposition}{Proposition}[section]
\newtheorem{lemma}{Lemma}[section]
\newtheorem{cor}{Corollary}[section]
\newtheorem{fact}{Fact}[section]
\newtheorem{remark}{Remark}[section]

\makeatletter
\@addtoreset{equation}{section}

\makeatother
\title
[Pfaffian of $F_4$]
{Pfaffian of Appell's hypergeometric system $F_4$ 
in terms of the intersection form of twisted cohomology groups
}
\author{Yoshiaki Goto}
\address[Goto]{
Department of Mathematics,
Graduate School of Science,
Kobe University,
Kobe 657-8501, Japan}
\email{y-goto@math.kobe-u.ac.jp}
\author{Jyoichi Kaneko}
\address[Kaneko]{
   Department of Mathematical Science,
   University of the Ryukyus,
   Nishihara, Okinawa, 903-0213, Japan
}
\email{kaneko@math.u-ryukyu.ac.jp}
\author{Keiji Matsumoto}
\address[Matsumoto]{
Department of Mathematics\\
Hokkaido University\\
Sapporo 060-0810, Japan
}
\email{matsu@math.sci.hokudai.ac.jp}

\keywords{Appell's hypergeometric system of 
 differential equations, Pfaffian system, 
Twisted cohomology group
}
\subjclass[2010]{33C65,  58A17.}
\date{\today}

\begin{document}

\begin{abstract}
We study a Pfaffian of the system of differential equations annihilating 
Appell's hypergeometric series $F_4(a,b,c;x)$ 
by twisted cohomology groups associated with 
integrals representing solutions to this system.
We simplify its connection matrix by the pull-back under 
a double cover of the complement of the singular locus. 
We express the simplified connection matrix 
in terms of the intersection form of the twisted cohomology groups.
\end{abstract}
\maketitle

\section{Introduction}
Appell's hypergeometric system $\cF_4(a,b,c)$  is generated by 
differential equations  (\ref{eq:DE}) annihilating 
Appell's hypergeometric series $F_4(a,b,c;x)$ of $x=(x_1,x_2)$ 
defined by  (\ref{eq:F4}), 
where $a,b$ and $c=(c_1,c_2)$ are parameters.
It is a holonomic system of rank $4$ with singular locus 
$$S=\{(x_1,x_2)\in \C^2\mid x_1x_2(x_1^2+x_2^2-2x_1x_2-2x_1-2x_2+1)=0\}
\cup L_\infty,$$
where $L_\infty$ is the line at infinity in the projective plane $\P^2$.

In this paper, we study a Pfaffian of the system 
$\cF_4(a,b,c)$ by a twisted cohomology group $H^2(\W^\bu(\C_x^2),\na)$ 
associated to integrals (\ref{eq:int-rep}) representing solutions 
to this system.
We regard the Pfaffian of $\cF_4(a,b,c)$ 
as a connection $\na_X$ of a vector bundle over $X=\P^2-S$ 
with fiber $H^2(\W^\bu(\C_x^2),\na)$ over $x\in X$.
We find a frame of this vector bundle and represent the connection 
as $\Xi=\Xi^1dx_1+\Xi^2dx_2$ in Theorem \ref{th:Pfaff} in terms of 
$4\times 4$-matrices $\Xi^1$ and $\Xi^2$.
The connection matrix satisfies the integrability condition 
$d\Xi =\Xi\wedge \Xi$ and has no apparent singularity,
but $\Xi_1$ and $\Xi_2$ are complicated and $d\Xi\ne O$. 
To make the connection matrix simple, we consider 
a double cover $Y$ of $X$  defined by a map 
$$\C^2\ni (y_1,y_2)\mapsto (x_1,x_2)=(y_1(1-y_2),y_2(1-y_2))\in \C^2.$$
It induces the pull-back bundle  and the pull-back connection $\na_Y$.
By changing a frame of the pull-back bundle, 
we express the pull-back connection $\na_Y$ 
by $\wh \Xi$ in terms of logarithmic $1$-forms of 
the variables $y_1$ and $y_2$ in Theorem \ref{th:conn-hat}. 
In particular, we have $d\wh \Xi=\wh \Xi\wedge \wh \Xi=O$.  
It is studied in \cite{Kat1} that the pull-back of $F_4(a,b,c;x)$ under 
a similar map satisfies a Pfaffian equation equivalent to ours.

There is the intersection pairing $\cI_c$ between $H^2(\W^\bu(\C_x^2),\na)$ 
and $H^2(\W^\bu(\C_x^2),\na^\vee)$, which is defined by 
the dual derivative $\na^\vee$ of $\na$. 
We have the dual vector bundle over $X$ and the dual connection $\na_X^\vee$
satisfying
$$d_x\cI_c(\f,\f')=\cI_c(\na_X\f,\f')+\cI_c(\f,\na_X^\vee\f'),$$
where $\f$ and $\f'$ are a section of the vector bundle 
and that of its dual, respectively, and 
$d_x$ is the exterior derivative with respect to $x_1$ and $x_2$.   
We can not find a frame of $H^2(\W^\bu(\C_x^2),\na)$ satisfying 
\begin{equation}
\label{eq:flat}
d_x\cI_c(\f_i,\f_j^\vee)=0
\end{equation}
for any $1\le i,j\le4$, where $\f_i\in H^2(\W^\bu(\C_x^2),\na)$ 
and $\f_j^\vee$ is the image of $\f_j$ under the natural map 
from $H^2(\W^\bu(\C_x^2),\na)$ to $ H^2(\W^\bu(\C_x^2),\na^\vee)$. 
The intersection pairing is defined between the pull-back bundles, and 
the frame of the pull-back bundle used in Theorem \ref{th:conn-hat} 
satisfies the relations (\ref{eq:flat}) for any $1\le i,j\le4$.
The existence of the frame enables us to represent the pull-back 
connection $\na_Y$  by the intersection form $\cI_c$ 
in Theorem \ref{th:exp-conn-int}.  
We remark that this expression is not given in terms of matrices.

The monodromy representation of $\cF_4(a,b,c)$ is 
initially given in \cite{K} by a 
twisted homology group associated to integrals (\ref{eq:int-rep}),  
and reconsidered in \cite{GM} by the intersection pairing.
In \cite{GM}, we characterize the circuit transform of a loop turning around 
each component of $S$ by the intersection form and a subspace of 
vanishing cycles as $x$ approaches the component.
In Lemma \ref{lem:roots-exp}, we characterize the coefficient matrix of 
the logarithmic $1$-form corresponding to a component of 
$\wt{S}$ in $\wh \Xi$   by the intersection form and a subspace of 
vanishing forms as $x$ approaches the component,
where $\wt{S}$ is the preimage of $S$ under the double cover. 

Appell's hypergeometric system $\cF_4(a,b,c)$ 
is generalized to Lauricella's hypergeometric system $\cF_C(a,b,c)$ 
of rank $2^m$ with $m$ variables. For this system, 
we have twisted (co)homology groups associated to integrals 
representing solutions. 
The monodromy representation of $\cF_C(a,b,c)$ is expressed 
in terms of intersection form between twisted homology groups in \cite{G}.
However, we have not seen that the system $\cF_C(a,b,c)$ 
admits a Pfaffian system with a simple expression of 
the pull-back connection under some covering map. 
We think that it is important to find a covering map and 
a frame of its pull-back bundle satisfying the relations 
(\ref{eq:flat}) for any $1\le i,j \le 2^m$.

For studies of Pfaffians of Lauricella's hypergeometric systems 
$\cF_A$ and $\cF_D$ of $m$-variables, it is easy to find such frames 
without considering covering maps. Their Pfaffians are expressed in terms of 
intersection forms between twisted cohomology groups associated to integrals 
representing solutions, refer to \cite{M1} and \cite{M2}.


\section{Appell's hypergeometric function $F_4$}
In \cite{ApKa}, Appell's hypergeometric series $F_4(a,b,c;x)$ 
of variables $x_1,x_2$ with parameters 
$a$, $b$, $c=(c_1,c_2)$ is defined by 
\begin{equation}
\label{eq:F4}
F_4(a,b,c_1,c_2;x_1,x_2)=\sum_{n_1,n_2=0}^\infty
\frac{(a,n_1+n_2)(b,n_1+n_2)}{(c_1,n_1)(c_2,n_2)n_1!n_2!}x_1^{n_1}x_2^{n_2},
\end{equation}
where $c_1,c_2\ne 0,-1,-2,\dots$ and 
$(a,k)=a(a+1)\cdots(a+k-1)=\G(a+k)/\G(a)$.
This series converges in the domain 
$$\D=\{x=(x_1,x_2)\in \C\mid \sqrt{|x_1|}+\sqrt{|x_2|}<1\},$$
and satisfies differential equations 
\begin{eqnarray}
\nonumber& &
\Big[x_1(1- x_1)\pa_1^2- x_2^2\pa_2^2- 2x_1x_2\pa_1\pa_2\\
\nonumber& &\ + \{c_1- (a+ b+ 1)x_1\}\pa_1- 
(a+ b+ 1)x_2\pa_2- ab\Big]f(x)=0,
\\
\label{eq:DE}
\\
\nonumber
& &
\Big[x_2(1- x_2)\pa_2^2- x_1^2\pa_1^2- 2x_1x_2\pa_1\pa_2\\ 
\nonumber
& &\ + \{c_2- (a+ b+ 1)x_2\}\pa_2- 
(a+ b+ 1)x_1\pa_1- ab\Big]f(x)=0.
\end{eqnarray}
The system generated by them is called Appell's hypergeometric 
system $\mathcal{F}_4(a,b,c_1,c_2)$ of differential equations.
This system is of rank $4$ with singular locus
$$S=\{(x_1,x_2)\in \C^2\mid x_1x_2 R(x)=0\}\cup L_\infty\subset \P^2,$$
where $L_\infty$ is the line at infinity and 
$$R(x)=x_1^2+x_2^2-2x_1x_2-2x_1-2x_2+1.$$

We use the following integral representation of solutions of 
$\mathcal{F}_4(a,b,c)$ 
\begin{equation}
\label{eq:int-rep}
\int_\De u(s,x) \frac{ds_1\wedge ds_2}{s_1s_2},
\end{equation}
where $$
u(s,x)=s_1^{\lm_1}s_2^{\lm_2}Q(s)^{\lm_3}L(s,x)^{\lm_4},
$$
$$Q=Q(s)=s_1s_2-s_1-s_2,\quad L=L(s,x)=1-s_1x_1-s_2x_2,$$
$$\lm_1=a-c_2+1,\ \lm_2=a-c_1+1,\  \lm_3=c_1+c_2-a-2,\ \lm_4=-b,$$ 
and a $2$-chain $\De$ loading a branch of $u(s,x)$ is a twisted $2$-cycle. 
Refer to \cite{AoKi} for twisted cycles and twisted homology groups. 
\begin{remark}
\label{rem:diff-para}
It is convenient for the study of a Pfaffian system of $\cF_4(a,b,c)$ 
to use the multi-valued function $u(s,x)$ in stead of 
$$u(t,x)=t_1^{1-c_1}t_2^{1-c_2}(1-t_1-t_2)^{c_1+c_2-a-1}
(1-\frac{x_1}{t_1}-\frac{x_2}{t_2})^{-b}
$$ 
used in \cite{GM}. We have $u(s,x)$  by the change of variables
$(t_1,t_2)=(1/s_1,1/s_2)$ and the replacement  $a\mapsto a+1$ 
for $u(t,x)$. 
\end{remark}

\section{Twisted cohomology group}
We regard the parameters $a,b$ and $c=(c_1,c_2)$ as indeterminants and 
we set 
$$
\begin{array}{cccc}
a_{00}=a,& a_{10}=a-c_1+1, &a_{01}=a-c_2+1, &a_{11}=a-c_1-c_2+2, \\
b_{00}=b,& b_{10}=b-c_1+1, &b_{01}=b-c_2+1, &b_{11}=b-c_1-c_2+2. \\
\end{array}
$$
We assume that 
\begin{equation}
\label{eq:non-integer}
a_{ij},b_{ij}\notin \Z\quad (i,j\in \F_2=\{0,1\})
\end{equation}
when we assign complex values to the parameters.
Recall that 
$$\lm_1=a-c_2+1,\ \lm_2=a-c_1+1,\ \lm_3=c_1+c_2-a-2,\ \lm_4=-b.$$ 
In this section, we regard vector spaces as defined over 
the rational function field $\C(\lm)=\C(\lm_1,\dots,\lm_4)=\C(a,b,c_1,c_2)$.

We set 
$$\frak X=\big\{(s,x)\in \C^2\times X\big|s_1s_2Q(s)L(s,x)\ne0\big\}
\subset (\P^1\times \P^1)\times \P^{2},$$
where  $X$ is the complement of the singular locus $S$ of $\cF_4(a,b,c)$ in 
$\P^2$. 
There is a natural projection 
$$\pr:\frak X\ni (s,x)\mapsto x\in X.$$ 
For any fixed $x\in X$, we have 
$$\C_x^2=\pr^{-1}(x)=\{s=(s_1,s_2)\in \C^2\mid s_1s_2Q(s)L(s,x)\ne0\}$$
and an inclusion map 
$$\imath_x:\C_x^2\ni s \mapsto (s,x)\in \frak{X}.$$ 


We denote the $\C(\lm)$-algebra of rational functions on $\P^2$ 
with poles only along $S$ by $\cO(X)$. 
We denote the vector space of rational $k$-forms on 
$(\P^1\times \P^1)\times \P^{2}$ 
with poles only along the complement of $\frak X$ 
by $\W^k(\frak X)$ and 
the subspace of $\W^{i+j}(\frak X)$
consisting elements which are $i$-forms with respect to the variables 
$s_1,s_2$ by $\W^{i,j}(\frak X)$.

We set 
\begin{eqnarray*}
\w&=&d_s\log(u(s,x))=\frac{\lm_1ds_1}{s_1}+\frac{\lm_2ds_2}{s_2}+
\frac{\lm_3d_sQ(s)}{Q(s)}+\frac{\lm_4d_sL(s,x)}{L(s,x)},\\
\w_X&=&d_x\log(u(s,x))=-\frac{\lm_4s_1dx_1}{L(s,x)}
-\frac{\lm_4s_2dx_2}{L(s,x)},
\end{eqnarray*}
where $d_s$ and $d_x$ are the exterior derivative with respect to 
$s_1,s_2$ and to $x_1,x_2$, respectively.
Note that $\w\in \W^{1,0}(\frak X)$ and 
$\w_X\in \W^{0,1}(\frak X)$.
By a twisted exterior derivation $\na=d_s+\w\wedge$ on $\frak X$,
we define quotient spaces
$$
H^k(\W^{\bu,0}(\frak X),\na)=\ker\left(\na:\W^{k,0}(\frak X)\to 
\W^{k+1,0}(\frak X)\right)
\big/\na\left(\W^{k-1,0}(\frak X)\right)$$
as $\cO(X)$-modules, 
where $k=0,1,2$ and we regard $\W^{-1,0}(\frak X)$ as 
the zero vector space.

For a fixed $x$, the inclusion map $\imath_x$ induces a natural map 
from $H^k(\W^{\bu,0}(\frak X),\na)$ to 
the rational twisted cohomology group
$$H^k(\W^\bu(\C_x^2),\na)=
\ker\left(\na:\W^k(\C_x^2)\to \W^{k+1}(\C_x^2)\right)
/\na\left(\W^{k-1}(\C_x^2)\right)$$
on $\C_x^2$ with respect to the twisted exterior derivative induced from 
$\na$. 
Here $\W^k(\C_x^2)$ is the vector space of rational $k$-forms 
with poles only along the complement of $\C_x^2$ in $\P^1\times \P^1$.
The structure of $H^k(\W^\bu(\C_x^2),\na)$ is known as follows.
\begin{fact}[\cite{AoKi},\cite{Cho}]
\label{fact:TCH}
\begin{itemize}
\item[$\mathrm{(i)}$] We have 
$$\dim H^k(\W^\bu(\C_x^2),\na)=\left\{
\begin{array}{ccl}
4 &\textrm{if}&k=2,\\
0 &\textrm{if}&k=0,1.
\end{array}
\right.
$$
\item[$\mathrm{(ii)}$]
There is a canonical isomorphism  
\begin{eqnarray*}
\jmath_x: H^2(\W^\bu(\C_x^2),\na) 
&\to& H^2(\cE_c^\bu(\C_x^2),\na)\\
&=&\ker(\na:\cE_c^{2}(\C_x^2)\to \cE_c^{3}(\C_x^2))
/\na(\cE_c^{1}(x)),
\end{eqnarray*}
where $\cE_c^k(\C_x^2)$ is the vector space of smooth $k$-forms 
with compact support in $\C_x^2$.
\end{itemize}
\end{fact}

We have a twisted exterior derivation $\na^\vee=d_s-\w\wedge$
for $-\w$ and 
\begin{eqnarray*}
H^2(\W^{\bu,0}(\frak X),\na^\vee)&=&
\W^{2,0}(\frak X)\big/\na^\vee(\W^{1,0}(\frak X)),\\
H^2(\W^\bu(\C_x^2),\na^\vee)&=&
\W^2(\C_x^2)/\na^\vee(\W^{1}(\C_x^2)).
\end{eqnarray*}
For any fixed $x\in X$, we define the intersection form between 
$H^2(\W^\bu(\C_x^2),\na)$
and 
$H^2(\W^\bu(\C_x^2),\na^\vee)$
by 
$$\cI_c(\f_x,\f_x')=\int_{\C_x^2}\jmath_x(\f_x)\wedge \f_x'\in \C(\lm),$$
where $\f_x,\f_x'\in \W^{2}(\C_x^2)$,  $\jmath_x$ is given in 
Fact \ref{fact:TCH}. 
This integral converges since $\jmath_x(\f_x)$ is a smooth $2$-form 
on $\C_x^2$ with compact support. 
It is bi-linear over $\C(\lm)$.

We take four elements $\f_1,\dots,\f_4$ of $H^2(\W^{\bu,0}(\frak X),\na)$: 
$$\begin{array}{ll}
\f_1=\dfrac{ds_{12}}{s_1s_2},&
\f_2=\dfrac{x_1ds_{12}}{s_2L(s,x)},\\[5mm]
\f_3=\dfrac{x_2ds_{12}}{s_1L(s,x)},& 
\f_4=\dfrac{ds_{12}}{Q(s)L(s,x)},
\end{array}
$$
where $ds_{12}=ds_1\wedge ds_2.$

\begin{proposition}
\label{prop:int.no}
For a fixed $x\in X$, the numbers 
$\cI_c(\imath_x^*(\f_i),\imath_x^*(\f_j))$ $(1\le i,j\le 4)$ are 
$(2\pi\sqrt{-1})^2C_{ij}$, where 
$$\begin{array}{ccl}
C_{11}&=&
\dfrac{1}{\lm_{123}}\left(
\dfrac{1}{\lm_{1}}+\dfrac{1}{\lm_{2}}\right)
+
\dfrac{1}{\lm_{134}^-}\left(
\dfrac{1}{\lm_{0}}+\dfrac{1}{\lm_{2}}\right)
+
\dfrac{1}{\lm_{234}^-}\left(
\dfrac{1}{\lm_{0}}+\dfrac{1}{\lm_{1}}\right)\\
&=&\dfrac{1}{a_{00}}\left(
\dfrac{1}{a_{01}}+\dfrac{1}{a_{10}}
\right)+
\dfrac{1}{b_{10}}\left(
\dfrac{1}{b_{11}}+\dfrac{1}{a_{10}}
\right)+
\dfrac{1}{b_{01}}\left(
\dfrac{1}{b_{11}}+\dfrac{1}{a_{01}}
\right),\\
C_{12}&=&
\dfrac{-1}{\lm_{134}^-}
\left(\dfrac{1}{\lm_{0}}+\dfrac{1}{\lm_{2}}\right)
=\dfrac{-1}{b_{10}}
\left(\dfrac{1}{b_{11}}+\dfrac{1}{a_{10}}\right)
,\\
C_{13}&=& 
\dfrac{-1}{\lm_{234}^-}
\left(\dfrac{1}{\lm_{0}}+\dfrac{1}{\lm_{1}}\right)
=\dfrac{-1}{b_{01}}\left(
\dfrac{1}{b_{11}}+\dfrac{1}{a_{01}}
\right)
,\\
C_{14}\!&\!=\!&0,\\
C_{22}\!&\!=\!& 
\left(\dfrac{1}{\lm_{0}}+\dfrac{1}{\lm_{2}}\right)
\left(\dfrac{1}{\lm_{4}}+\dfrac{1}{\lm_{134}^-}\right)
=\left(\dfrac{1}{b_{11}}+\dfrac{1}{a_{10}}\right)
\left(\dfrac{-1}{b_{00}}+\dfrac{1}{b_{10}}\right)
,\\
C_{23}&=& 
\dfrac{-1}{\lm_0\lm_{4}}=\dfrac{1}{b_{11}b_{00}}
,\\
C_{24}\!&\!=\!&0,\\
C_{33}\!&\!=\!& 
\left(\dfrac{1}{\lm_{0}}+\dfrac{1}{\lm_1}\right)
\left(\dfrac{1}{\lm_{4}}+\dfrac{1}{\lm_{234}^-}\right)
=\left(\dfrac{1}{b_{11}}+\dfrac{1}{a_{01}}\right)
\left(\dfrac{-1}{b_{00}}+\dfrac{1}{b_{01}}\right)
,\\
C_{34}\!&\!=\!&0,\\
C_{44}\!&\!=\!&\dfrac{2}{\lm_3\lm_4R(x)}=\dfrac{2}{a_{11}b_{00}R(x)},\\
C_{ji}\!&\!=\!& C_{ij} \quad \textrm{for }\ i<j,
\end{array}
$$
$$
\begin{array}{ll}
\lm_0=-(\lm_1+\lm_2+2\lm_3+\lm_4)=b_{11},
&\lm_{123}=\lm_1+\lm_2+\lm_3=a_{00},\\
 \lm_{134}^-=-(\lm_1+\lm_3+\lm_4)=b_{10},
&\lm_{234}^-=-(\lm_2+\lm_3+\lm_4)=b_{01}.
\end{array}
$$
The matrix $C=(C_{ij})_{i,j}$ is symmetric and its determinant is 
$$\det(C)=\dfrac{4\lm_3}{\lm_0\lm_1\lm_2\lm_4^3
\lm_{123}\lm_{134}^-\lm_{234}^-R(x)}=
\left(\prod_{i,j=0,1}\frac{1}{a_{ij}b_{ij}}\right)
\frac{4a_{11}^2}{b_{00}^2R(x)}.$$
\end{proposition}
\proof 
For a fixed $x$, we blow up $\P^1\times \P^1(\supset \C_x^2)$ 
at two points $(0,0)$ and $(\infty,\infty)$. 
We tabulate the residue of the pull back of $\w$ to this space 
in Table \ref{tab:res}, where $E_0$ and $E_\infty$ are exceptional divisors 
coming from the blow up of $\P^1\times \P^1$ at 
the points $(0,0)$ and $(\infty,\infty)$, respectively.
\begin{table}[hbt]
$$
\begin{array}{|c|c|r|c|}
\hline
\textrm{divisor} & \lm\textrm{'s} &a_{ij},b_{ij} & a,b,c_1,c_2\\
\hline 
E_\infty& \lm_0 & b_{11} & b-c_1-c_2+2 \\ 
s_1=0& \lm_1 & a_{01} & a-c_2+1 \\ 
s_2=0& \lm_2 & a_{10} & a-c_1+1 \\ 
Q(s)=0& \lm_3 &-a_{11} &-a+c_1+c_2-2 \\ 
L(s,x)=0& \lm_4 &-b_{00} & -b \\ 
E_{0}& \lm_{123}& a_{00} & a \\ 
s_1=\infty& \lm_{134}^-& b_{10} & b-c_1+1 \\ 
s_2=\infty& \lm_{234}^-& b_{01} & b-c_2+1 \\ 
\hline
\end{array}
$$ 
\caption{
Residues for components of the pole divisor of $\w$
}
\label{tab:res}
\end{table}
\begin{figure}[hbt]
  \centering
\includegraphics[width=10cm]{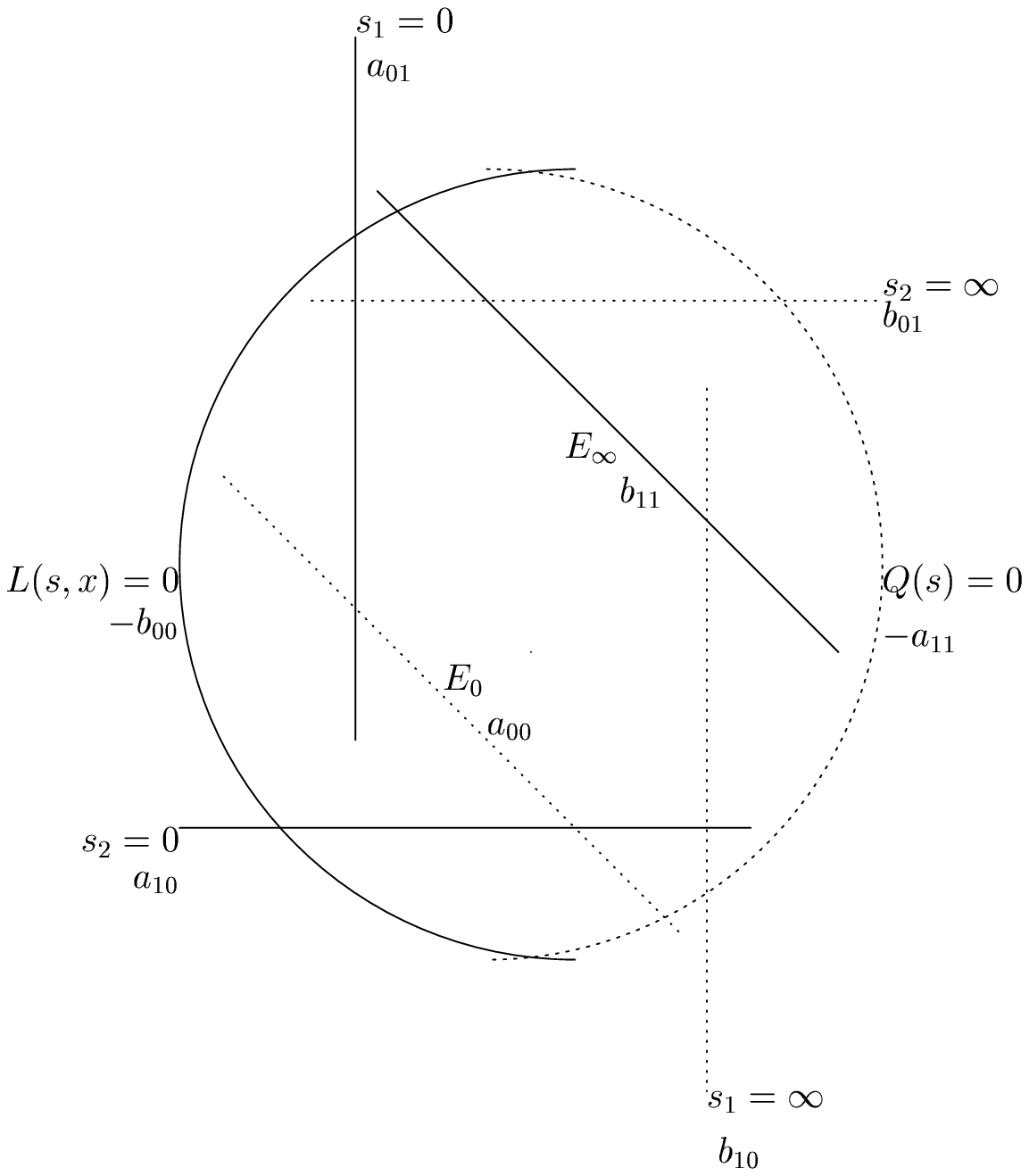}  
  \caption{Pole divisor of $\w$}
  \label{fig:pole-div}
\end{figure}
By using these data, follow the proof of Theorem 5.1 in \cite{GM}. 
\Qed

Note that the matrix $C$ is well-defined and 
$\det(C)\ne0$ for any $x\in X$ under our assumption.
The natural map from $H^2(\W^{\bu,0}(\frak X),\na)$ to 
$H^2(\W^\bu(\C_x^2),\na)$ is surjective by Fact \ref{fact:TCH}.

\begin{cor}
\label{cor:frame}
We can regard the $\cO(X)$-modules $H^2(\W^{\bu,0}(\frak X),\na)$  
and $H^2(\W^{\bu,0}(\frak X),\na^\vee)$ 
as vector bundles 
$$
\bigcup_{x\in X}H^2(\W^\bu(\C_x^2),\na),\quad 
\bigcup_{x\in X}H^2(\W^\bu(\C_x^2),\na^\vee),
$$
over $X$.  
The classes of $\f_1,\dots,\f_4$ form frames of 
$H^2(\W^{\bu,0}(\frak X),\na)$ and $H^2(\W^{\bu,0}(\frak X),\na^\vee)$. 
\end{cor}

\proof
We have only to prove that 
the natural map from $H^2(\W^{\bu,0}(\frak X),\na)$ to 
$H^2(\W^\bu(\C_x^2),\na)$ is injective. 
We show that 
if 
an element $\f\in \W^{2,0}(\frak X)$ satisfying 
$\imath_x^*(\f)=0$ as an element of $H^2(\W^\bu(\C_x^2),\na)$ 
for any fixed $x\in X$ then 
$\f$ belongs to $\na(\W^{1,0}(\frak X))$.
There exists $\psi_x\in \W^1(\C_x^2)$ such that 
$\na\psi_x=\imath_x^*(\f)$ for any $x$. 
Since this differential equation is of variables $s_1,s_2$ 
with parameters $x_1,x_2$, $\psi_x$ can be globally extended to $\psi$.  
Hence we have $\psi\in \W^{1,0}(\frak X)$ such that 
$\na(\psi)=\f$. 
\Qed


By Proposition \ref{prop:int.no} and  Corollary \ref{cor:frame}, 
the intersection form $\cI_c$ is regarded as 
a map from $H^2(\W^{\bu,0}(\frak X),\na)\times 
H^2(\W^{\bu,0}(\frak X),\na^\vee)$ to $\cO(X)$.

\section{Pfaffian system of $\cF_4(a,b,c)$}
 
For any $\f\in H^2(\W^{\bu,0}(\frak X),\na)$, we have
\begin{equation}
\label{eq:d2nax}
d_x\int_\De u(s,x)\f
=\int_\De (d_xu(s,x)\wedge \f+u(s,x)d_x\f)=\int_\De u(s,x)(\na_X \f),
\end{equation}
where a $2$-chain $\De$ loading a branch of $u(s,x)$ is a twisted $2$-cycle. 
%
Thus the exterior derivative $d_x$ on $X$ induces 
the connection $\na_X=d_x+\w_X\wedge$
$$\na_X:
H^2(\W^{\bu,0}(\frak X),\na)\to 
H^2(\W^{\bu,1}(\frak X),\na)
=\W^{2,1}(\fX)/\na(\W^{1,1}(\fX)).
$$ 
By considering $1/u(s,x)$ instead of $u(s,x)$, we also have the connection 
$\na_X^\vee =d_x-\w_X\wedge$
$$\na_X^\vee:
H^2(\W^{\bu,0}(\frak X),\na^\vee)\to 
H^2(\W^{\bu,1}(\frak X),\na^\vee)
=\W^{2,1}(\fX)/\na^\vee(\W^{1,1}(\fX)).
$$ 
\begin{proposition}
\label{prop:int-str}
The connection $\na_X$ is compatible with the intersection form $\cI_c$, i.e., 
they satisfy 
$$d_x\cI_c(\f,\f')=\cI_c(\na_X \f,\f')+\cI_c(\f,\na_X^\vee \f').$$
\end{proposition} 
\proof
It is enough to show this identity in any small simply connected domain $U$ 
in $X$. 
We have 
\begin{eqnarray*}
& &d_x\cI_c(\f,\f')
=d_x\int_{\C_x^2} u(s,x)\jmath_x(\f)\wedge\frac{\f'}{u(s,x)}\\
&=&\int_{\C_x^2} d_x[u(s,x)\jmath_x(\f)]\wedge \frac{\f'}{u(s,x)}
+\int_{\C_x^2} u(s,x)\jmath_x(\f)\wedge d_x\frac{\f'}{u(s,x)}\\
&=&\int_{\C_x^2} \na_X(\jmath_x(\f))\wedge \f'
+\int_{\C_x^2} \jmath_x(\f)\wedge \na_X^\vee \f',
\end{eqnarray*}
since $\jmath_x(\f)$ is with compact support for any point $x\in U$. 
By following the proof of Lemma 7.2 in \cite{M1}, we can show that 
$\na_X(\jmath_x(\f))$ is cohomologue to $\jmath_x(\na_X(\f))$.
Hence we have 
$$
\int_{\C_x^2} \na_X(\jmath_x(\f))\wedge \f'
=\int_{\C_x^2} \jmath_x(\na_X(\f))\wedge \f'
=\cI_c(\na_X\f,\f'),
$$
which completes the proof. 
\Qed

Since $\na_X \f\in \W^{2,1}(\fX)$, there exist 
$\Xi_{i}^1$ and $\Xi_{i}^2$ in $\cO(X)$ such that 
$$\na_X \f=
 dx_1\wedge \sum_{i=1}^4 \Xi^1_i\f_i+dx_2\wedge\sum_{i=1}^4 \Xi^2_i\f_i.
$$
By calculating $\Xi_i^1$ and $\Xi_i^2$ for $\f=\f_1,\dots,\f_4$, 
we represent the connection $\na_X$ as
$$\na_X\tr(\f_1,\dots,\f_4)=\Xi\wedge \tr(\f_1,\dots,\f_4),\quad 
\Xi=dx_1\Xi^1+dx_2 \Xi^2,$$
where $\Xi^1$ and $\Xi^2$ are $4\times 4$-matrices over the algebra $\cO(X)$.
By (\ref{eq:d2nax}),  the vector valued function 
$$F(x)=\tr\left(\int_\De u(s,x)\f_1,\dots,\int_\De u(s,x)\f_4\right)
$$
satisfies a system of differential equations 
$$ 
d_xF(x)= \Xi \wedge F(x).
$$
Let us determine $\Xi^1$ and $\Xi^2$. A straightforward calculation implies 
the following. 


\begin{lemma}
\label{lem:diff}
\begin{eqnarray*}
\na_X(\f_1)&=&dx_1\wedge \frac{-\lm_4ds_{12}}{s_2L}
+dx_2\wedge\frac{-\lm_4ds_{12}}{s_1L} \\
&=&\frac{dx_1}{x_1}\wedge(-\lm_4)\f_2+\frac{dx_2}{x_2}\wedge(-\lm_4)\f_3,\\
\na_X(\f_2)&=&dx_1 \wedge \frac{(1-\lm_4s_1x_1-s_2x_2)ds_{12}}{s_2L^2}
+dx_2\wedge \frac{(1-\lm_4)x_1 ds_{12}}{L^2},\\
\na_X(\f_3)&=&dx_1\wedge \frac{(1-\lm_4)x_2 ds_{12}}{L^2}
+dx_2\wedge \frac{(1-s_1x_1-\lm_4s_2x_2)ds_{12}}{s_1L^2},\\
\na_X(\f_4)&=&
 d x_1\wedge \frac{(1-\lm_4)s_1ds_{12}}{QL^2}
+d x_2\wedge \frac{(1-\lm_4)s_2ds_{12}}{QL^2}.
\end{eqnarray*}
\end{lemma}
To obtain $\Xi^1$ and $\Xi^2$, we express 
$$\frac{ds_{12}}{L^2},\quad
\frac{s_1ds_{12}}{s_2L^2},\quad 
\frac{s_2ds_{12}}{s_1L^2},\quad 
\frac{ds_{12}}{s_2L^2},\quad 
\frac{ds_{12}}{s_1L^2},\quad 
\frac{s_1ds_{12}}{QL^2},\quad 
\frac{s_2ds_{12}}{QL^2},$$
in terms of $\f_1,\dots,\f_4$.

\begin{lemma}
\label{lem:L2}
As elements of $H^2(\W^{\bu,0}(\fX),\na)$, we have
\begin{eqnarray*}
\frac{2(\lm_4-1)x_1x_2}{L^2}ds_{12}&=&
(\lm_1\!+\!\lm_2\!+\!\lm_3)\f_1+(\lm_2\!-\!\lm_4)\f_2\\
& &+(\lm_1\!-\!\lm_4)\f_3+\lm_3(1\!-\! x_1\!-\! x_2)\f_4,\\
\frac{(\lm_4-1)x_1s_1}{s_2L^2}ds_{12}
&=&\frac{\lm_1+\lm_3+1}{x_1}\f_2+\lm_3\f_4,\\
\frac{(\lm_4-1)x_2s_2}{s_1L^2}ds_{12}
&=&\frac{\lm_2+\lm_3+1}{x_2}\f_3+\lm_3\f_4.
\end{eqnarray*}

\end{lemma}
\proof
Straightforward calculations imply 
\begin{eqnarray*}
& &\!\!
\na\left(\frac{x_1ds_1-x_2ds_2}{L}+\frac{ds_1}{s_1}-\frac{ds_2}{s_2}\right)\\
&=&\!\!\frac{2(\lm_4-1)x_1x_2ds_{12}}{L^2}\\
& &\!\!-[(\lm_1\!+\!\lm_2\!+\!\lm_3)\f_1+(\lm_2\!-\!\lm_4)\f_2
+(\lm_1\!-\!\lm_4)\f_3
+\lm_3(1\!-\! x_1\!-\! x_2)\f_4],
\end{eqnarray*}
\begin{eqnarray*}
& &\na\left(-\frac{s_1ds_2}{s_2L}\right)
=\frac{(\lm_4-1)x_1s_1}{s_2L^2}ds_{12}
-\big[\frac{\lm_1+\lm_3+1}{x_1}\f_2+\lm_3\f_4\big],\\
& &\na\left(\frac{s_2ds_1}{s_1L}\right)
=\frac{(\lm_4-1)x_2s_2}{s_1L^2}ds_{12}
-\big[\frac{\lm_2+\lm_3+1}{x_2}\f_3+\lm_3\f_4\big],
\end{eqnarray*}
which shows this lemma.
\Qed

\begin{lemma}
\label{lem:1/sL2}
As elements of $H^2(\W^{\bu,0}(\fX),\na)$, we have
\begin{eqnarray*}
 \frac{2(\lm_4-1)x_2}{s_1L^2}ds_{12}&=&
(\lm_1\!+\!\lm_2\!+\!\lm_3)\f_1+(\lm_2\!-\!\lm_4)\f_2\\
& &+(\lm_1\!+\!2\lm_2\!+\!2\lm_3\!+\!\lm_4)\f_3+\lm_3(1\!-\! x_1\!+\! x_2)\f_4,
\\
\frac{2(\lm_4-1)x_1}{s_2L^2}ds_{12}&=&
(\lm_1\!+\!\lm_2\!+\!\lm_3)\f_1+(2\lm_1\!+\!\lm_2\!+\!2\lm_3\!+\!\lm_4)\f_2\\
& &+(\lm_1\!-\!\lm_4)\f_3+\lm_3(1\!+\! x_1\! -\! x_2)\f_4.
\end{eqnarray*}
\end{lemma}
\proof
Note that 
\begin{eqnarray*}
& &\left(\frac{1}{s_1L^2}-\frac{x_1}{L^2}-\frac{s_2x_2}{s_1L^2}\right)ds_{12}
=\frac{ds_{12}}{s_1L}=\frac{\f_3}{x_2},\\
& &\left(\frac{1}{s_2L^2}-\frac{x_2}{L^2}-\frac{s_1x_1}{s_2L^2}\right)ds_{12}
=\frac{ds_{12}}{s_2L}=\frac{\f_2}{x_1}.
\end{eqnarray*}
Use Lemma \ref{lem:L2}.
\Qed

\begin{lemma}
\label{lem:QL2}
\begin{eqnarray*}
\frac{s_1ds_{12}}{QL^2}
&=&\frac{1\!-\! x_1\!+\! x_2}{R(x)}\f_4
+\frac{(\lm_1+\lm_2+2\lm_3+\lm_4)(1+x_1-x_2)x_2}{\lm_3R(x)}
\frac{ds_{12}}{L^2}\\
& &
\!\!
-\frac{2\lm_1x_2}{\lm_3R(x)}\frac{ds_{12}}{s_1L^2}
-\frac{\lm_2(1-x_1-x_2)}{\lm_3R(x)}
\frac{ds_{12}}{s_2L^2},\\
\frac{s_1ds_{12}}{QL^2}
&=&\frac{1\!+\! x_1\!-\! x_2}{R(x)}\f_4
+\frac{(\lm_1+\lm_2+2\lm_3+\lm_4)(1-x_1+x_2)x_1}{\lm_3R(x)}\frac{ds_{12}}{L^2}
\\
& &
-\frac{\lm_1(1-x_1-x_2)}{\lm_3R(x)}\frac{ds_{12}}{s_1L^2}
-\frac{2\lm_2x_1}{\lm_3R(x)}\frac{ds_{12}}{s_2L^2}.
\end{eqnarray*}
\end{lemma}

\proof
Set
$$\eta_0=\frac{ds_{12}}{QL^2},\quad \eta_1=\frac{s_1ds_{12}}{QL^2},\quad 
\eta_2=\frac{s_2ds_{12}}{QL^2}.$$
There is a relation
$$\eta_0-x_1\eta_1-x_2\eta_2=\frac{(1-s_1x_1-s_2x_2)ds_{12}}{QL^2}=
\f_4$$
among them. We have                  
\begin{eqnarray*}
& &\na\left(\frac{x_1ds_1+x_2ds_2}{L^2}\right)\\
&=&\lm_3(x_1-x_2)\eta_0-\lm_3x_1\eta_1+\lm_3x_2\eta_2
+(\frac{\lm_1x_2}{s_1L^2}-\frac{\lm_2x_1}{s_2L^2})ds_{12},\\
& &\na\left(\frac{ds_1-ds_2}{L^2}+\frac{(x_1+x_2)(-s_2ds_1+s_1ds_2)}{L^2}
\right)\\
&=&2\lm_3\eta_0+\lm_3(x_1+x_2-1)\eta_1+\lm_3(x_1+x_2-1)\eta_2\\
& &-\left(\frac{\lm_1}{s_1L^2}+\frac{\lm_2}{s_2L^2}-
\frac{(\lm_1+\lm_2+2\lm_3+\lm_4)(x_1+x_2)}{L^2}\right)ds_{12},
\end{eqnarray*}
which are zero as elements of  $H^2(\W^{\bu,0}(\frak X),\na)$.
By regarding these relations as linear equations of variables $\eta_0$,
$\eta_1$, $\eta_2$, solve them. 
\Qed

\begin{remark} The form $\eta_0={ds_{12}}/{QL^2}$ 
in the proof of Lemma \ref{lem:QL2} is expressed as
\begin{eqnarray*}
\eta_0
&=&
\frac{1-x_1-x_2}{R(x)}\f_4+
\frac{2(\lm_1+\lm_2+2\lm_3+\lm_4)x_1x_2}{\lm_3 R(x)}
\frac{ds_{12}}{L^2}\\
& &-\frac{\lm_1x_2(1+x_1-x_2)}{\lm_3R(x)}
\frac{ds_{12}}{s_1L^2}
   -\frac{\lm_2x_1(1-x_1+x_2)}{\lm_3R(x)}
\frac{ds_{12}}{s_2L^2}.
\end{eqnarray*}

\end{remark}

Lemmas \ref{lem:diff}, 
\ref{lem:L2}, \ref{lem:1/sL2} and \ref{lem:QL2} yield
the representation of $\na_X$ 
with respect to the frame $\tr(\f_1,\dots,\f_4)$.  
\begin{theorem}
\label{th:Pfaff}
With respect to the frame $\tr(\f_1,\dots,\f_4)$, 
the connection $\na_X$ is represented as
$$\na_X \tr(\f_1,\dots,\f_4)=\Xi\wedge \tr(\f_1,\dots,\f_4),
$$
where $\Xi=\Xi^1dx_1+\Xi^2dx_2$ and 
$$\Xi^1=\begin{pmatrix}
0&-\dfrac{\lm_4}{x_1}&0&0\\
0&-\dfrac{\lm_1+\lm_3}{x_1}&0&-\lm_3\\
-\dfrac{\lm_1+\lm_2+\lm_3}{2x_1}&
-\dfrac{\lm_2-\lm_4}{2x_1}&
-\dfrac{\lm_1-\lm_4}{2x_1}&
-\dfrac{\lm_3 (1-x_1-x_2)}{2x_1}\\[3mm]
\Xi^1_{4,1}&\Xi^1_{4,2}&\Xi^1_{4,3}&\Xi^1_{4,4}\\
\end{pmatrix},
$$
$$
\Xi^2=\begin{pmatrix}
0&0&-\dfrac{\lm_4}{x_2}&0\\
-\dfrac{\lm_1+\lm_2+\lm_3}{2x_2}&
-\dfrac{\lm_2-\lm_4}{2x_2}&
-\dfrac{\lm_1-\lm_4}{2x_2}&
-\dfrac{\lm_3 (1-x_1-x_2)}{2x_2}\\
0&0&-\dfrac{\lm_2+\lm_3}{x_2}&-\lm_3\\[3mm]
\Xi^2_{4,1}&\Xi^2_{4,2}&\Xi^2_{4,3}&\Xi^2_{4,4}
\end{pmatrix},
$$
\begin{eqnarray*}
\Xi^1_{4,1}
&=&\frac{\lm_1+\lm_2+\lm_3}{\lm_3}\left(\frac{\lm_1-\lm_2}{R(x)}-
\frac{(\lm_1+2\lm_3+\lm_4)(1+x_1-x_2)}{2x_1R(x)}\right),\\
\Xi^1_{4,2}&=&-\frac{(\lm_2-\lm_4)(\lm_2+2\lm_3+\lm_4)}{\lm_3R(x)}\\
& &+\frac{\lm_1\lm_2+\lm_1\lm_4+\lm_2\lm_4+2\lm_3\lm_4+\lm_4^2}{2\lm_3}
\frac{1-x_1-x_2}{x_1R(x)},\\
\Xi^1_{4,3}
&=&\frac{\lm_1\lm_2+\lm_1\lm_4+\lm_2\lm_4+2\lm_3\lm_4+\lm_4^2}{\lm_3R(x)}\\
& &-\frac{(\lm_1-\lm_4)(\lm_1+\lm_4+2\lm_3)}{2\lm_3}
\frac{1-x_1-x_2}{x_1R(x)},\\
\Xi^1_{4,4}
&=&-(\lm_1-2\lm_3-3\lm_4+2)\frac{x_1}{2R(x)}+
(\lm_1-\lm_4+1)\frac{1+x_2}{R(x)}\\
& &-(\lm_1+2\lm_3+\lm_4)\frac{(x_2-1)^2}{2x_1R(x)},\\
\end{eqnarray*}
\begin{eqnarray*}
\Xi^2_{4,1}
&=&\frac{\lm_1+\lm_2+\lm_3}{\lm_3}
\left(\frac{\lm_2-\lm_1}{R(x)}
-\frac{(\lm_2+2\lm_3+\lm_4)(1-x_1+x_2)}{2x_2R(x)}\right),\\
\Xi^2_{4,2}&=&
\frac{\lm_1\lm_2+\lm_1\lm_4+\lm_2\lm_4+2\lm_3\lm_4+\lm_4^2}{\lm_3R(x)}\\
& &-\frac{(\lm_2-\lm_4)(\lm_2+2\lm_3+\lm_4)}{2\lm_3}\frac{1-x_1-x_2}{x_2R(x)},\\
\Xi^2_{4,3}&=&-\frac{(\lm_1-\lm_4)(\lm_1+\lm_4+2\lm_3)}{\lm_3R(x)}\\
& &
+\frac{\lm_1\lm_2+\lm_1\lm_4+\lm_2\lm_4+2\lm_3\lm_4+\lm_4^2}{2\lm_3}
\frac{1-x_1-x_2}{x_2R(x)},\\
\Xi^2_{4,4}&=&
-(\lm_2-2\lm_3-3\lm_4+2)\frac{x_2}{2 R(x)}+
(\lm_2-\lm_4+1)\frac{x_1+1}{R(x)}\\
& &-(\lm_2+2\lm_3+\lm_4)\frac{(x_1-1)^2}{2x_2 R(x)}.
\end{eqnarray*}
\end{theorem}
\begin{remark}
Though the connection matrix $\Xi$ is not closed, 
it satisfies the integrability condition , i.e., 
\begin{eqnarray*}
& &\Xi\wedge \Xi=(\Xi^1\Xi^2-\Xi^2\Xi^1)dx_1\wedge dx_2\\
&=&d_x \Xi =
\left(-\frac{\pa}{\pa x_2}\Xi^1+\frac{\pa}{\pa x_1}\Xi^2\right)
dx_1\wedge dx_2\ne O.
\end{eqnarray*}
\end{remark}
\begin{remark}
\label{rem:cohomolog}
We give some expressions of elements of $H^2(\W^{\bu,0}(\frak X),\na)$
in terms of $\f_1,\dots,\f_4$:
\begin{eqnarray*}
\na(-\frac{ds_1}{s_1}+\frac{ds_2}{s_2})
&=&\frac{\lm_3ds_{12}}{Q}-[-(\lm_1+\lm_2+\lm_3)\f_1+\lm_4\f_2+\lm_4\f_3],\\
\na(\frac{ds_2}{s_2})
&=&\frac{\lm_3ds_{12}}{s_1Q}-[-(\lm_1+\lm_3)\f_1+\lm_4\f_2],\\
\na(-\frac{ds_1}{s_1})
&=&\frac{\lm_3ds_{12}}{s_2Q}-[-(\lm_2+\lm_3)\f_1+\lm_4\f_3].
\end{eqnarray*}
\end{remark}

\begin{remark}
\label{rem:Pfaff}
We define a function $f_1(x)$ by $\int_\De u(s,x)\f_1$ 
for a twisted cycle $\De$ loading a branch of $u(s,x)$. 
This function satisfies 
$$x_1\frac{\pa f_1(x)}{\pa x_1}=-\lm_4\int_\De u(s,x)\f_2,\quad 
x_2\frac{\pa f_1(x)}{\pa x_2}=-\lm_4\int_\De u(s,x)\f_3,$$
\begin{eqnarray*}
x_1x_2\frac{\pa^2f_1(x)}{\pa x_1\pa x_2}
&=&
\frac{\lm_4}{2}\Big[
\lm_{123}\!\int_\De u(s,x)\f_1\!+\!(\lm_2\!-\!\lm_4)\!\int_\De u(s,x)\f_2\\
& &\hspace{4mm}+\!(\lm_1\!-\!\lm_4)\!\int_\De u(s,x)\f_3
\!+\!\lm_3(1\!-\! x_1\!-\! x_2)\!\int_\De u(s,x)\f_4\Big].
\end{eqnarray*}
Thus the vector-valued function 
$$F_\pa(x)=\tr\left(f_1(x),x_1\frac{\pa f_1(x)}{\pa x_1},
x_2\frac{\pa f_1(x)}{\pa x_2},x_1x_2\frac{\pa^2 f_1(x)}{\pa x_1\pa x_2}\right)
$$
satisfies 
$$d_x F_\pa(x)=(G_\pa \Xi G_\pa^{-1}+d_xG_\pa\;G_\pa^{-1})F_\pa(x),$$
where 
$$G_\pa=\left(\begin{array}{cccc}
1 & 0 & 0 & 0 \\
0 &-\lm_4 & 0 & 0 \\
0 & 0 & -\lm_4  & 0 \\
\frac{\lm_{123}\lm_4}{2} & \frac{(\lm_2-\lm_4)\lm_4}{2} 
& \frac{(\lm_1-\lm_4)\lm_4}{2}  & \frac{\lm_3\lm_4(1-x_1-x_2)}{2}\\
  \end{array}\right),
$$  
$$G_\pa^{-1}=\left(\begin{array}{cccc}
1 & 0 & 0 & 0 \\
0 &-\lm_4^{-1} & 0 & 0 \\
0 & 0 & -\lm_4^{-1}  & 0 \\
\frac{-\lm_{123}}{\lm_3(1-x_1-x_2)} &\frac{\lm_2-\lm_4}{\lm_3\lm_4(1-x_1-x_2)} 
&\frac{\lm_1-\lm_4}{\lm_3\lm_4(1-x_1-x_2)}  
& \frac{2}{\lm_3\lm_4(1-x_1-x_2)}\\
  \end{array}\right).
$$  
Note that the matrix $G_\pa$ does not belongs to $GL_4(\cO(X))$.
\end{remark}

\section{Connection matrix in terms of intersection form}

We define an affine variety 
$$\cX=\{(x_1,x_2,x_3)\in \C^3\mid x_3^2= R(x)\},$$ 
which is regarded  as the double cover of $\C^2$ branching along 
the divisor $R(x_1,x_2)=0$ by  the projection 
$$\pro:\cX\ni (x_1,x_2,x_3)\to (x_1,x_2)\in \C^2.$$
Note that 
$$\left(\frac{1+x_1-x_2-x_3}{2}\right)\left(\frac{1+x_1-x_2+x_3}{2}\right)=x_1,
$$
$$
\left(\frac{1-x_1+x_2-x_3}{2}\right)\left(\frac{1-x_1+x_2+x_3}{2}\right)=x_2,
$$
for $(x_1,x_2,x_3)\in \cX$. Thus the preimages of lines 
$x_1=0$ and $x_2=0$ in $\C^2$ under the projection $\pro$ are expressed as
 equations
$$\frac{1+x_1-x_2-x_3}{2}=0, \quad \frac{1-x_1+x_2-x_3}{2}=1,$$
and
$$\frac{1+x_1-x_2-x_3}{2}=1, \quad \frac{1-x_1+x_2-x_3}{2}=0,$$
in $\cX$ respectively.
By a map 
$$\cX\ni (x_1,x_2,x_3)\mapsto  
(y_1,y_2)=\left(\frac{1\!+\! x_1\!-\! x_2\!-\! x_3}{2},
\frac{1\!-\! x_1\!+\! x_2\!-\! x_3}{2}\right)\in 
\C^2,$$
$\cX$ is bi-holomorphic to $\C^2$; its inverse is expressed as 
$$\C^2\ni (y_1,y_2)\mapsto (x_1,x_2,x_3)=
(y_1(1\!-\! y_2),(1\!-\! y_1)y_2, 1\!-\! y_1\!-\! y_2)\in \cX.$$
Though $(x_1,x_2)$ are not valid as local coordinates
on the set 
$$\left\{(x_1,x_2,x_3)\in \cX\Big|
\frac{\pa}{\pa x_3}(x_3^2-R(x))=2x_3=0\right\},$$
we can use
$(y_1,y_2)=\left(\dfrac{1\!+\! x_1\!-\! x_2\!-\!x_3}{2},
\dfrac{1\!-\! x_1\!+\! x_2\!-\! x_3}{2}\right)$
as a global coordinates system on $\cX$.
The covering transformation 
$$\rho:(x_1,x_2,x_3)\mapsto (x_1,x_2,-x_3)$$
of $\pro:\cX\to X$ is represented as 
$$(y_1,y_2)\mapsto (1-y_2,1-y_1).$$
The ramification locus of $\pro$ is 
$$\{(y_1,y_2)\in \C^2\mid y_1+y_2=1\}$$
and the preimage of the singular locus of $\cF_4(a,b,c)$ in $\C^2$ 
under the projection $\pro$ is 
$$\{(y_1,y_2)\in \C^2\mid [y_1(1-y_2)]\cdot[y_2(1-y_1)]\cdot(1-y_1-y_2)=0\}.$$
We set 
\begin{eqnarray*}
Y&=&\{y=(y_1,y_2)\in\C^2\mid y_1(1-y_1)y_2(1-y_2)(1-y_1-y_2)\ne 0\}
\subset\P^2,\\
\wt{S}&=&\P^2-Y.
\end{eqnarray*}
Note that $\wt{S}$ and  $Y$ 
are invariant as sets under the action of $\rho$.
By the vector bundle $H^2(\W^{\bu,0}(\frak X),\na)$ over $X$ and 
the projection $pr:Y\to X$, we have the vector bundle 
$$pr^*H^2(\W^{\bu,0}(\frak X),\na)=
\{(y,\f)\in Y\times H^2(\W^{\bu,0}(\frak X),\na)\mid pr(y)=p(\f)\}.
$$
Its frame $\tr(\wt\f_1,\dots,\wt\f_4)$ is given by 
the pull-backs 
$$\wt\f_i=pr^*(\f_i)\quad (i=1,\dots,4)$$
of $\f_1,\dots,\f_4$ in $H^2(\W^{\bu,0}(\frak X),\na)$
under the projection $pr$.
The connection $\na_Y=pr^*\na_X$ of $pr^*H^2(\W^{\bu,0}(\frak X),\na)$ 
for the frame $\tr(\wt\f_1,\dots,\wt\f_4)$
is defined by the pull-back of $\Xi$ under the projection $pr$.

\begin{cor}
\label{cor:pullback}
With respect to the frame $\tr(\wt\f_1,\dots,\wt\f_4)$, 
the connection $\na_{Y}$
is expressed as 
$$
\wt\Xi \wedge \tr(\wt\f_1,\dots,\wt\f_4),\quad \wt\Xi=\wt\Xi^1 dy_1+\wt\Xi^2 dy_2,$$
where 
$$
\wt\Xi^1=\left(
\begin {array}{cccc} 0&{\frac {-\lm_4}{ y_1}}
&
{\frac {-\lm_4}{ y_1-1}}&0\\[2mm] 
-{\frac {\lm_{123}}{2({ y_1}-1)}}
&
{\frac {-\lm_2+\lm_4}{2(y_1-1)}}-{\frac {\lm_1+\lm_3}{ y_1}}
&{\frac {-\lm_1+\lm_4}{2(y_1-1)}}
&-\frac{ \lm_3}{2}-{\frac {{ \lm_3}{ y_2}}{2({ y_1}-1)}}
\\[2mm]
-{\frac {\lm_{123}}{2 y_1}}
&{\frac {-\lm_2+\lm_4}{2 y_1}}
&{\frac {-\lm_1+\lm_4}{2y_1}}-{\frac {\lm_2+\lm_3}{ y_1-1}}
&\frac{\lm_3}{2}+{\frac { \lm_3(y_2-1)}{2y_1}}\\[2mm]
\wt\Xi^1_{41}&\wt\Xi^1_{42}&\wt\Xi^1_{43}&\wt\Xi^1_{44}
\end {array}\right),
$$
$$
\wt\Xi^2=\left(
\begin {array}{cccc} 
0&{\frac {-\lm_4}{y_2-1}}&{\frac {-\lm_4}{y_2}}&0
\\[2mm] 
-{\frac {\lm_{123}}{2y_2}}
&{\frac {-\lm_2+\lm_4}{2y_2}}-{\frac {\lm_1+\lm_3}{y_2-1}}
&{\frac {-\lm_1+\lm_4}{2y_2}}
&\frac{\lm_3}{2}+{\frac {\lm_3(y_1 -1 )}2y_2}
\\[2mm] 
-{\frac {\lm_{123}}{2(y_2-1)}}
&{\frac {-\lm_2+\lm_4}{2(y_2-1)}}
&{\frac {-\lm_1+\lm_4}{2(y_2-1)}}-{\frac {\lm_2+\lm_3}{{y_2}}}
&-\frac{\lm_3}{2}-{\frac {\lm_3y_1}{2(y_2-1)}}
\\[2mm]
\wt\Xi^2_{41}&\wt\Xi^2_{42}&\wt\Xi^2_{43}&\wt\Xi^2_{44}
\end {array}
\right),
$$

\begin{eqnarray*}
\wt\Xi^1_{41}&=&
{\frac { ( \lm_3-\lm_{134}^- ) 
  \lm_{123}  }{2\lm_3y_1
 (y_2 -1 ) }}
-\frac { ( \lm_3-\lm_{234}^- ) 
  \lm_{123} }
{2\lm_3(y_1 -1 ) y_2}
\\
& &+{\frac { (\lm_3-\lm_{234}^- )\lm_{123}}{2\lm_3y_2 (y_1+y_2-1 )}}
-{\frac { (2\lm_3+{ \lm_4}+{ \lm_1} )
 \lm_{123}}{2\lm_3 (y_2-1 ) (y_1+y_2-1 )}},\\[2mm]
\wt\Xi^1_{42}&=&
{\frac {(\lm_2+\lm_4)(\lm_1+\lm_4)}{2\lm_3(y_2-1)(y_1+y_2-1)}}
-{\frac {(\lm_2+\lm_4)(\lm_1+\lm_4)}{2\lm_3(y_2-1)y_1}}
+{\frac {{\lm_4}}{(y_2-1)(y_1+y_2-1)}}\\
& &
-{\frac {\lm_4}{y_1(y_2-1)}}
-{\frac {(\lm_2-\lm_4)(\lm_3-\lm_{234}^-)}{2\lm_3y_2(y_1-1)}}
+{\frac {(\lm_2-\lm_4)(\lm_3-\lm_{234}^-)}{2\lm_3y_2(y_1+y_2-1)}},
\\[2mm]
\wt\Xi^1_{43}&=&
%
{\frac {(\lm_2+\lm_4)(\lm_1+\lm_4)}{2\lm_3y_2(y_1-1)}}
-{\frac {(\lm_2+\lm_4)(\lm_1+\lm_4)}{2\lm_3y_2(y_1+y_2-1)}}
-{\frac {\lm_4}{y_2(y_1+y_2-1)}}\\
& &+{\frac {\lm_4}{y_2(y_1-1)}}
-{\frac {(\lm_1-\lm_4)(\lm_3-\lm_{134}^-)}{2\lm_3(y_2-1)(y_1+y_2-1)}}
+{\frac {(\lm_1-\lm_4)(\lm_3-\lm_{134}^-)}{2\lm_3(y_2-1)y_1}},
\\[2mm]
\wt\Xi^1_{44}&=&{\frac {2\lm_4+2\lm_3-1}{y_1+y_2-1}}
-{\frac {\lm_3-\lm_{234}^-}{2(y_1-1)}}
-{\frac {\lm_3-\lm_{134}^-}{2y_1}},
\end{eqnarray*}

\begin{eqnarray*}
\wt\Xi^2_{41}&=&-{\frac { (2 \lm_3+\lm_4+\lm_1 ) 
\lm_{123}}{2\lm_3 y_1 (y_2-1)}}
+{\frac { (2 \lm_3+\lm_4+\lm_1 )\lm_{123}}{2\lm_3 y_1  (y_1+y_2-1 )}}
\\
& &
-{\frac { (\lm_2+2 \lm_3+\lm_4 ) \lm_{123} }
{2\lm_3  (y_1-1) (y_1+y_2-1 )}}
+{\frac { (\lm_2+2 \lm_3+\lm_4 ) \lm_{123}}
{2\lm_3y_2  (y_1-1)}},
\\[2mm]
\wt\Xi^2_{42}&=&\frac{ (\lm_2+\lm_4 ) (\lm_1+\lm_4 )}
{2\lm_3 (y_2-1 )y_1}
-\frac{ (\lm_2+\lm_4 ) (\lm_1+\lm_4 )}{2\lm_3y_1 (y_1+y_2-1 )}
-\frac{\lm_4}{y_1 (y_1+y_2-1 )}\\
& &+\frac{\lm_4}{ (y_2-1 )y_1}
-{\frac { (\lm_2-\lm_4 ) (\lm_3-\lm_{234}^- )}
{2\lm_3 (y_1-1 ) (y_1+y_2-1 )}}
+{\frac { (\lm_2-\lm_4 ) (\lm_3-\lm_{234}^- )}
{2\lm_3y_2 (y_1-1 )}},
\\[2mm]
\wt\Xi^2_{43}&=&
-\frac{ (\lm_2+\lm_4 ) (\lm_1+\lm_4 )}{2\lm_3y_2 (y_1-1 )}
+\frac{ (\lm_2+\lm_4 ) (\lm_1+\lm_4 )}{2\lm_3 (y_1-1 ) (y_1+y_2-1 )}
+\frac{\lm_4}{ (y_1-1 ) (y_1+y_2-1 )}
\\
& &
-\frac{\lm_4}{y_2 (y_1-1 )}
-{\frac { (\lm_1-\lm_4 ) (\lm_3-\lm_{134}^- )}
{2\lm_3 (y_2-1 )y_1}}
+{\frac { (\lm_1-\lm_4 ) (\lm_3-\lm_{134}^- )}
{2\lm_3y_1 (y_1+y_2-1 )}}
,
\\[2mm]
\wt\Xi^2_{44}&=&
-{\frac {\lm_3-\lm_{134}^-}{2(y_2-1)}}
+{\frac {2\lm_3+2\lm_4-1}{y_1+y_2-1}}
-{\frac {\lm_3-\lm_{234}^-}{2y_2}}.
\end{eqnarray*}
\end{cor}

\proof
By using
$$\begin{array}{ll}
\pro^*(x_1)=y_1(1-y_2), &\pro^*(x_1)=(1-y_1)y_2,\\
\pro^*(dx_1)=(1-y_2)dy_1-y_1dy_2,&
\pro^*(dx_2)=-y_2dy_1+(1-y_1)dy_2,\\
\pro^*(R(x))=(1-y_1-y_2)^2,
\end{array}
$$
we have only to calculate the pull back of $\Xi$ under the projection $\pro$. 
\Qed

Let $\cO(Y)$ be the $\C(\lm)$-algebra of rational functions on 
$\P^2$ with poles only along the complement of $Y$.
We change the frame  $\tr(\wt\f_1,\dots,\wt\f_4)$ 
of the vector bundle $pr^*H^2(\W^{\bu,0}(\frak X),\na)$ 
over $Y$ 
as 
$$\tr(\wh\f_1,\dots,\wh\f_4)=G \tr(\wt\f_1,\dots,\wt\f_4),$$
$$
G=\left(\begin{array}{cccc}
1& & & \\
 &1& & \\
 & &1& \\
 & & &1-y_1-y_2
  \end{array}
\right)\in GL_4(\cO(Y)).  
$$
The  connection matrix $\wh{\Xi}$ of $pr^*H^2(\W^{\bu,0}(\frak X),\na)$ 
with respect to the frame $\tr(\wh\f_1,\dots,\wh\f_4)$  is given by 
the gauge transformation 
$$G\;\wt \Xi\; G^{-1}+d_yG\; G^{-1}$$
of $\wt \Xi$.  Straightforward calculations 
imply the following.

\begin{theorem}
\label{th:conn-hat}
The connection matrix $\wh\Xi$ is 
$$\wh\Xi^1\frac{dy_1}{y_1}+\wh\Xi^2\frac{dy_2}{y_2}
+I_{3,1}\wh\Xi^2I_{3,1}^{-1}\frac{dy_1}{y_1-1}
+I_{3,1}\wh\Xi^1I_{3,1}^{-1}\frac{dy_2}{y_2-1}
+\wh\Xi^3\frac{dy_1+dy_2}{y_1+y_2-1},$$
where   
\begin{eqnarray*}
\wh\Xi^1&=&\left(\begin{array}{cccc} 
0&-\lm_4&0&0
\\
0&-{\lm_1}-\lm_3&0&0
\\
\frac{-\lm_{123}}{2}&\frac{-\lm_2+\lm_4}{2}&\frac{-\lm_1+\lm_4}{2}
&\frac{-\lm_3}{2}
\\
{\frac {(\lm_{134}^- -\lm_3)\lm_{123}}{2{\lm_3}}}
&\lm_4\!+\!{\frac {(\lm_2+\lm_4 ) (\lm_1+\lm_4 )}{2\lm_3}}
&
{\frac {(\lm_1-\lm_4)(\lm_{134}^- -\lm_3)}{2\lm_3}}
&\frac{\lm_{134}^- -{\lm_3}}{2}\end {array}\right),
\\
\wh\Xi^2&=&\left(\begin {array}{cccc} 
0&0&-\lm_4&0\\
\frac{-\lm_{123}}{2}&\frac{-\lm_2+\lm_4}{2}
&\frac{-\lm_1+\lm_4}{2}&\frac{-\lm_3}{2}
\\\noalign{\medskip}0&0&-\lm_3-\lm_2&0\\
\noalign{\medskip}
{\frac {(\lm_{234}^- -\lm_3)\lm_{123}}{2\lm_3}}
&
{\frac{(\lm_2-\lm_4)(\lm_{234}^- -\lm_3)}{2\lm_3}}
&\lm_4\!+\!{\frac {(\lm_2+\lm_4)(\lm_1+\lm_4)}{2\lm_3}}
&\frac{\lm_{234}^- -\lm_3}{2}
\end {array}\right ),
\\
\wh\Xi^3&=&
\left(\begin {array}{cccc} 
0&0&0&0\\
0&0&0&0\\
0&0&0&0\\
0&0&0&2(\lm_3+\lm_4)
\end {array}\right ),
\qquad I_{3,1}=
\left(\begin {array}{cccc} 
1&0&0&0\\
0&1&0&0\\
0&0&1&0\\
0&0&0&-1
\end {array}\right ).
\end{eqnarray*}
The connection matrix $\wh\Xi$  satisfies 
$$d_y\wh\Xi =\wh\Xi \wedge \wh\Xi=O.$$

\begin{remark}
\label{rem:closed}
The gauge transformation by $G$ changes 
the non-closed connection matrix $\wt\Xi$ 
into the closed connection matrix $\wh\Xi$.
\end{remark}

\end{theorem}

We express the connection $\na_Y$ in terms of intersection form $\cI_c$.
Since the map $pr:Y\to X$ is locally isomorphic, we have 
$$\cI_c(\wt\f_i,\wt\f_j)=\cI_c(\f_i,\f_j),\quad (1\le i,j\le 4).$$
We give the intersection numbers for the frame $\tr(\wh\f_1,\dots,\wh\f_4)$.

\begin{cor}
\label{cor:int.no}
The intersection numbers $\cI_c(\wh\f_i,\wh\f_j)$ $(1\le i,j\le 4)$ 
are $(2\pi\sqrt{-1})^2\wh{C}_{ij}$,  where 
$$\wh{C}_{ij}=\left\{\begin{array}{cl}
C_{ij},  &\textrm{if }1\le i,j\le 3,\\[2mm]
\dfrac{2}{\lm_3\lm_4}=\dfrac{2}{a_{11}b_{00}},\quad&
\textrm{if } (i,j)=(4,4),\\[2mm]
0, & \textrm{otherwise},\\
\end{array}
\right.
$$
$C_{ij}$ are in Proposition \ref{prop:int.no}. 
\end{cor}
\proof It is clear by the transformation $G$, 
$$\pro^*(R(x))=(1-y_1-y_2)^2,$$
and Proposition \ref{prop:int.no}. 
\Qed

Let $\wh C$ be the matrix $(\wh C_{ij})_{1\le i,j\le 4}$. 
This matrix is symmetric and $\wh C^\vee=\wh C$, 
where $\wh C^\vee$ is given by the replacement $\lm_i\to -\lm_i$ 
$(i=1,\dots,4)$ 
for every entry of $\wh C$.

\begin{lemma}
\label{lem:orth}
The connection matrix $\wh \Xi$ satisfies 
$$\wh\Xi^\vee=-\wh\Xi,\quad \wh\Xi \; \wh C+\wh C\;^t\wh\Xi^\vee=O,$$
where $\wh\Xi^\vee$ is given by the replacement $\lm_i\to -\lm_i$ 
for $\wh\Xi$.
\end{lemma}
\proof We can easily check this lemma by Theorem  
\ref{th:conn-hat} and Corollary \ref{cor:int.no}.
The second equality is also obtained by 
Proposition \ref{prop:int-str} from $d_y\wh C=O$. 
\Qed

\begin{lemma}
\label{lem:eigenspace}
\begin{itemize}
\item[(i)] 
The eigenvalues of $\wh\Xi^1$ are $0$ and $-(\lm_1+\lm_3)=1-c_1$ and 
each of the eigenspaces is $2$-dimensional.
Its $(1-c_1)$-eigenspace is spanned by the row vectors
$$e_2=(0,1,0,0),\quad 
(\frac{\lm_{123}}{\lm_3},\frac{\lm_2-\lm_4}{\lm_3},
\frac{\lm_1-\lm_4}{\lm_3},1).$$
\item[(ii)] 
The eigenvalues of $\wh\Xi^2$ are $0$ and $-(\lm_2+\lm_3)=1-c_2$ and 
each of the eigenspaces is $2$-dimensional.
Its $(1-c_2)$-eigenspace is spanned by the row vectors
$$e_3=(0,0,1,0),\quad 
(\frac{\lm_{123}}{\lm_3},\frac{\lm_2-\lm_4}{\lm_3},
\frac{\lm_1-\lm_4}{\lm_3},1).$$
\item[(iii)] 
The eigenvalues of $\wh\Xi^3$ are $0$ and $2(\lm_3+\lm_4)=2(c_1+c_2-a-b-2)$ 
and the $0$-eigenspace of $\wh\Xi^3$ is $3$-dimensional and 
the $2(\lm_3+\lm_4)$-eigenspace of $\wh\Xi^3$ is $1$-dimensional. 
Its $2(c_1+c_2-a-b-2)$-eigenvector is $e_4=(0,0,0,1).$
\end{itemize}
\end{lemma}
\proof  
We can easily check this lemma by Theorem  
\ref{th:conn-hat}.
\Qed

We set 
$$e_5=(\frac{\lm_{123}}{\lm_3},\frac{\lm_2-\lm_4}{\lm_3},
\frac{\lm_1-\lm_4}{\lm_3},1),\quad 
e_6=(\frac{\lm_{123}}{\lm_3},\frac{\lm_2-\lm_4}{\lm_3},
\frac{\lm_1-\lm_4}{\lm_3},-1),$$
and 
$$\wh\f_5=\frac{\lm_{123}}{\lm_3}\wh\f_1
+\frac{\lm_2-\lm_4}{\lm_3}\wh\f_2
+\frac{\lm_1-\lm_4}{\lm_3}\wh\f_3+\wh\f_4,$$ 
$$\wh\f_6=\frac{\lm_{123}}{\lm_3}\wh\f_1
+\frac{\lm_2-\lm_4}{\lm_3}\wh\f_2
+\frac{\lm_1-\lm_4}{\lm_3}\wh\f_3-\wh\f_4,$$ 
corresponding to the vectors $e_5$ and $e_6$, respectively.

\begin{lemma}
\label{lem:vanishi}
\begin{itemize}
\item[(i)]
The forms $\wh\f_2$ and $\wh\f_5$ vanish as $y_1\to 0$. 
The forms $\wh\f_2$ and $\wh\f_6$ vanish as $y_2\to 1$. 
\item[(ii)]
The forms $\wh\f_3$ and $\wh\f_5$ vanish as $y_2\to 0$. 
The forms $\wh\f_3$ and $\wh\f_6$ vanish as $y_1\to 1$. 
\item[(iii)]
The form $\wh\f_4$ vanishes as $y_1+y_2\to 1$. 
\end{itemize}
\end{lemma}

\proof
(i) Since $\wh\f_2=y_1(1-y_2)ds_{12}/(s_2L(s,y))$, it vanishes
as $y_1\to 0$ and as $y_2\to 1$.
We consider the image of $2(\lm_4-1)x_1x_2ds_{12}/L^2$ 
under the map $\pro^*$. By Lemma \ref{lem:L2}, 
it is 
$$\lm_{123}\wh\f_1+(\lm_2-\lm_4)\wh\f_2
+(\lm_1-\lm_4)\wh\f_3+\frac{\lm_3(1-y_1-y_2+2y_1y_2)}{1-y_1-y_2}\wh\f_4.
$$
It is clear that this element vanishes and 
its last term converges to $\lm_3\wh\f_4$ (resp. $-\lm_3\wh\f_4$)
as $y_1\to 0$ (resp. $y_2\to 1$). 
Thus $\lm_3\wh\f_5$ vanishes as $y_1\to 0$ and 
$\lm_3\wh\f_6$ vanishes as $y_2\to 1$.

\medskip\noindent
(ii) Similarly we can show the statements.

\medskip\noindent
(iii) Since 
$$\wh\f_4=\frac{(1-y_1-y_2)ds_{12}}{Q(s)L(s,y)},$$ 
it vanishes as $y_1+y_2\to 1$.
\Qed

We set 
$$\wh{C}_1=
\left(\begin{array}{cc}
\cI_c(\wh\f_2,\wh\f_2)&\cI_c(\wh\f_2,\wh\f_5)\\
\cI_c(\wh\f_5,\wh\f_2)&\cI_c(\wh\f_5,\wh\f_5)
\end{array}\right)
=\left(\begin{array}{cc}
\frac{(\lm_3-\lm_{134}^-)(\lm_1+\lm_3)}{\lm_0\lm_2\lm_4\lm_{134}^-}
&\frac{-2(\lm_1+\lm_3)}{\lm_0\lm_3\lm_4} \\
\frac{-2(\lm_1+\lm_3)}{\lm_0\lm_3\lm_4} 
&\frac{-4(\lm_2+\lm_3)(\lm_1+\lm_3)}{\lm_0\lm_3^2\lm_4}
\end{array}\right),
$$ 
$$\wh{C}_2=
\left(\begin{array}{cc}
\cI_c(\wh\f_3,\wh\f_3)&\cI_c(\wh\f_3,\wh\f_5)\\
\cI_c(\wh\f_5,\wh\f_3)&\cI_c(\wh\f_5,\wh\f_5)
\end{array}\right)
=\left(\begin{array}{cc}
\frac{(\lm_3-\lm_{234}^-)(\lm_2+\lm_3)}{\lm_0\lm_1\lm_4\lm_{234}^-}
&\frac{-2(\lm_2+\lm_3)}{\lm_0\lm_3\lm_4} \\
\frac{-2(\lm_2+\lm_3)}{\lm_0\lm_3\lm_4} 
&\frac{-4(\lm_2+\lm_3)(\lm_1+\lm_3)}{\lm_0\lm_3^2\lm_4}
\end{array}\right).
$$ 
We have
$$\left(\begin{array}{cc}
\cI_c(\wh\f_2,\wh\f_2)&\cI_c(\wh\f_2,\wh\f_6)\\
\cI_c(\wh\f_6,\wh\f_2)&\cI_c(\wh\f_6,\wh\f_6)
\end{array}\right)=\wh{C}_1,\ 
\left(\begin{array}{cc}
\cI_c(\wh\f_3,\wh\f_3)&\cI_c(\wh\f_3,\wh\f_6)\\
\cI_c(\wh\f_6,\wh\f_3)&\cI_c(\wh\f_6,\wh\f_6)
\end{array}\right)=\wh{C}_2.
$$

\begin{lemma}
\label{lem:roots-exp}
The matrix $\wh\Xi^i$ $(i=1,2,3)$ is expressed in terms of $\wh{C}$ and 
its eigenvectors with non-zero eigenvalue as 
\begin{eqnarray*}
\wh\Xi^1&=&-(\lm_1+\lm_3)\wh{C}(\tr e_2,\tr e_5)(\wh{C}_1)^{-1}
\left(\begin{array}{c} e_2\\ e_5  \end{array}\right),\\
\wh\Xi^2&=&-(\lm_2+\lm_3)\wh{C}(\tr e_3,\tr e_5)(\wh{C}_2)^{-1}
\left(\begin{array}{c} e_3\\ e_5  \end{array}\right),\\
\wh\Xi^3&=&2(\lm_3+\lm_4)\wh{C}\tr e_4(\wh{C}_{44})^{-1}e_4.
\end{eqnarray*}
\end{lemma}

\proof
We claim that 
$$v \wh{C} \tr w^\vee=0$$
for any eigenvector $v$ of $\wh \Xi^i$ with non-zero eigenvalue $\a$ 
and for any eigenvector $w$ of $\wh \Xi^i$ with eigenvalue $0$. 
Lemma \ref{lem:orth} implies 
$$\wh \Xi^i \wh C=-\wh C \tr (\wh\Xi^i)^\vee .$$
Since $\a\ne0$ and 
\begin{eqnarray*}
\a (v \wh C \tr w^\vee)
&=&(v\wh \Xi^i) \wh C \tr w^\vee 
=v(\wh\Xi^i \wh C) \tr w^\vee=-v (\wh C \tr (\wh\Xi^i)^\vee ) \tr w^\vee\\
&=&-v\wh C \tr( w \wh\Xi^i)^\vee =0,
\end{eqnarray*}
we have this claim.
Lemma \ref{lem:eigenspace} together with this claim 
gives the expression of the matrix $\wh\Xi^i$ in this lemma. 
\Qed

\begin{theorem}
\label{th:exp-conn-int}
The connection 
$\na_Y$ of $pr^*H^2(\W^{\bu,0}(\frak X),\na)$ 
is expressed in terms of the intersection form $\cI_c$ as 
\begin{eqnarray*}
 & &\na_{Y}(\wh\f)\\
&=&\hspace{3mm} \frac{dy_1}{y_1}\wedge 
(1-c_1)(\cI_c(\wh\f,\wh\f_2),\cI_c(\wh\f,\wh\f_5))(\wh C_1)^{-1}
\left(\begin{array}{c}
\wh\f_2\\ \wh\f_5
\end{array}\right)\\
&&+\frac{dy_2}{y_2}\wedge 
(1-c_2)(\cI_c(\wh\f,\wh\f_3),\cI_c(\wh\f,\wh\f_5))(\wh C_2)^{-1}
\left(\begin{array}{c}
\wh\f_3\\ \wh\f_5
\end{array}\right)\\
&&+\frac{dy_1}{y_1-1}\wedge 
(1-c_2)(\cI_c(\wh\f,\wh\f_3),\cI_c(\wh\f,\wh\f_6))(\wh C_2)^{-1}
\left(\begin{array}{c}
\wh\f_3\\ \wh\f_6
\end{array}\right)\\
&&+\frac{dy_2}{y_2-1}\wedge 
(1-c_1)(\cI_c(\wh\f,\wh\f_2),\cI_c(\wh\f,\wh\f_6))(\wh C_1)^{-1}
\left(\begin{array}{c}
\wh\f_2\\ \wh\f_6
\end{array}\right)\\
&&+\frac{dy_1+dy_2}{y_1+y_2-1}\wedge 
2(c_1+c_2-a-b-2)\cI_c(\wh\f,\wh\f_4)
(\wh C_{44})^{-1}\wh\f_4.
\end{eqnarray*}
\end{theorem}
\proof
Note that the linear transformation
$$\wh\f\mapsto (1-c_1)(\cI_c(\wh\f,\wh\f_2),\cI_c(\wh\f,\wh\f_5))(\wh C_1)^{-1}
\left(\begin{array}{c}
\wh\f_2\\ \wh\f_5
\end{array}\right)$$
is represented matrix $\wh\Xi^1$ with respect to the frame 
$\tr(\wh\f_1,\dots,\wh\f_4)$.
The eigen space of $I_{3,1}\wh\Xi^1I_{3,1}^{-1}$ 
with non-zero eigenvalue is spanned by $e_2$ and $e_6$, 
which correspond to $\wh\f_2$ and $\wh\f_6$, respectively. 
Thus the linear transformation 
$$\wh\f\mapsto (1-c_1)(\cI_c(\wh\f,\wh\f_2),\cI_c(\wh\f,\wh\f_6))(\wh C_1)^{-1}
\left(\begin{array}{c}
\wh\f_2\\ \wh\f_6
\end{array}\right)$$
is represented by $I_{3,1}\wh\Xi^1I_{3,1}^{-1}$ with respect to 
$\tr(\wh\f_1,\dots,\wh\f_4)$.
Similarly, we have representation matrices $\wh\Xi^2$ and 
$I_{3,1}\wh\Xi^2I_{3,1}^{-1}$.
The linear transformation 
$$2(c_1+c_2-a-b-2)\cI_c(\wh\f,\wh\f_4)
(\wh C_{44})^{-1}\wh\f_4$$ 
is represented by $\wh\Xi^3$ with respect to $\tr(\wh\f_1,\dots,\wh\f_4)$.
Theorem \ref{th:conn-hat} and Lemma \ref{lem:roots-exp} 
yield this theorem.
\Qed

\end{document}